\documentclass[a4paper,11pt]{amsart}
\date{\today}
\input epsf.sty
\usepackage{graphicx}
\usepackage{verbatim,enumerate}
\usepackage{amsthm}
\usepackage{boxedminipage}
\usepackage[matrix,arrow]{xy}
\usepackage{fancyheadings}
\usepackage{epic}

\theoremstyle{plain}
 \newtheorem{theorem}{Theorem}[section]
 \newtheorem*{theorem*}{Theorem}
 \newtheorem{proposition}[theorem]{Proposition}

\theoremstyle{definition}
 
 \newtheorem{conjecture}[theorem]{Conjecture}
\theoremstyle{remark}
 \newtheorem{remark}[theorem]{Remark}
 
\numberwithin{equation}{section}

\newcommand{\SO}{\mathrm{SO}}
\newcommand{\SU}{\mathrm{SU}}

\newcommand{\ad}{\mathrm{ad}}

\newcommand\Exp{\ensuremath{\mathrm{ Exp}}}

\title[On the Funk Transform]{On the Funk Transform on Compact Symmetric Spaces}

\author{Sebastian Klein}
\address[Klein]{%
Institut f\"ur Mathematik,
Universit\"at Mannheim,
68131 Mannheim,
Germany%
}
\email{s.klein@math.uni-mannheim.de}

\author{Gudlaugur Thorbergsson}
\address[Thorbergsson]{%
 Mathematisches Institut,
Universit\"at zu K\"oln,
Weyertal 86--90,
50931 K\"oln,
Germany%
}
\email{gthorbergsson@mi.uni-koeln.de}

\author{L\'aszl\'o Verh\'oczki}
\address[Verh\'oczki]{%
   E\"otv\"os University,
   Department of Geometry,
   P\'azm\'any P\'eter s. 1/C,
   1117 Budapest,
   Hungary
}
\email{verhol@cs.elte.hu}


\begin{document}
\begin{abstract} 
We prove that a function on an irreducible compact symmetric space $M$ which is not a sphere
is determined by its integrals over the shortest closed geodesics in $M$. 
We also prove a support theorem for the Funk transform on rank one symmetric spaces
which are not spheres.
\end{abstract}

\thanks{The first author's work was partially supported by a fellowship within the 
Postdoc-Programme of the German Academic Exchange Service (DAAD).
The second author was supported in part by the DFG-Schwerpunktprogramm {\em Globale
Differentialgeometrie.} Some of his work  was done during a visit 
at the Department of Geometry of the E\"otv\"os University in Budapest. He would like to thank 
this institution for the hospitality.
The third author was partially supported by the Hungarian
Scientific Research Fund OTKA K72537.} 

\maketitle
\section{Introduction}

Let $M$ be an irreducible simply connected compact symmetric space and denote by $\Xi$ the space
of the shortest closed geodesics on $M$. We associate to a smooth function $f$ on $M$ the function 
$\hat f$ on $\Xi$ defined by setting 
$$
\hat f(\gamma)=\int_\gamma f(\gamma(s))\,ds,
$$ 
where $s$ is the arc length parameter of $\gamma$. Following Helgason in \cite{He5}, we call the
linear map that sends $f$ to $\hat f$ the {\it Funk transform for $M$.}

The purpose of this paper is to prove the following two theorems.

\begin{theorem} \label{theorem}
The Funk transform for $M$ is injective if $M$ is not a sphere.
\end{theorem}

\begin{theorem} \label{support} Assume that $M$ has rank one and
is not a sphere. Let $\bar B_r(p)$ be a closed ball of radius $r>0$ around a point $p$ in $M$ 
where $r$ is smaller than the diameter of $M$. Let $f$ be a function on $M$ such that 
$\hat f(\gamma)=0$ for every $\gamma$ in $\Xi$ which does not meet $\bar B_r(p)$. 
Then $f(q)=0$ for all $q\not\in \bar B_r(p)$.
\end{theorem}

We recall that the simply connected rank one symmetric  spaces are precisely 
the spheres and the projective spaces over
the complex numbers, the quaternions and the octonions. 

The Funk transform of an odd function on a sphere $S^n$ clearly vanishes.
The following result, that we 
will  use in the proof of our theorems, was proved by Paul Funk 
in his doctoral dissertation and published in \cite{Fu} in 1913.

\begin{theorem}[Funk \cite{Fu}] \label{Funk}
Let $f$ be a function on a sphere $S^n$ with dimension $n\geq 2$
for which $\hat f=0$. Then $f$ is odd.
\end{theorem}

Other proofs of Theorem \ref{Funk} can be found in Appendix A of \cite{Gu} 
and in Chapter III, \S1B, of \cite{He4}.

The injectivity of the Funk transform is known for the simply connected projective spaces; see
\cite{He4}, p.~117. There it is based on Helgason's inversion formula for the antipodal Radon Transform which is difficult to prove; see Section 4 of Chapter I in \cite{He3} (or \cite{He6} where a 
substantially simplified proof based on ideas of Rouvi\`ere can be found). If one is only interested in
the injectivity of the Funk transform for these spaces, a simple proof can be given; see
 Remark \ref{simple}. 
 
In \cite{Gr},  Grinberg proved Theorem \ref{theorem} for the compact groups,  the complex and quaternionic Grassmannians, and $\SU(n)/\SO(n)$.

 Theorem \ref{support}  generalizes Theorem 3.4 in \cite{He5}. 
 
 There is a very interesting explicit inversion formula for the Funk transform
of functions with support in sufficiently small balls in  compact symmetric spaces; see \cite{He5}, Corollary 3.3.

\section{Helgason Spheres} \label{sec:helgspheres}

Let $M=G/K$ be an irreducible simply connected compact symmetric space.  If the metric of $M$ is
given by the negative of the Killing form, then the maximum of the sectional curvature on $M$ is
$\Vert\delta\Vert^2$, where $\delta$ is a highest restricted root.
We will normalize the metric on
$M$ in such a way that the maximum of the sectional curvature is equal to one. 
Then the injectivity radius satisfies $i(M)=\pi$.

The following
theorem is proved in \cite{He1}; see also \cite{He2}, Chapter VII, \S11.

\begin{theorem} \label{helgspheres}
 {\rm (1)} The shortest closed geodesics in $M$ have length $2\pi$.

{\rm (2)} Let $\gamma_1, \gamma_2:[0,1]\to M$ be two shortest closed geodesics. Then there is an
element $g$ in $G$ such that $g\circ\gamma_1=\gamma_2$.

{\rm (3)} A shortest closed geodesic in $M$ is contained in a totally geodesic sphere of constant curvature
one. The maximal dimension of such a sphere is $m(\delta)+1$, where $m(\delta)$ denotes the
multiplicity of a highest restricted root $\delta$. Any two such maximal spheres are conjugate 
under $G$.

\end{theorem}

The maximal totally geodesic spheres  in Theorem \ref{helgspheres} (3) are called {\it Helgason spheres}.
The projective lines in the simply connected projective spaces are Helgason spheres. They are used to 
give a simple proof of the injectivity of the Funk transform on projective spaces in Remark \ref{simple}.

It follows immediately from Theorem \ref{helgspheres} (2) that the space $\Xi$
of shortest closed geodesics in $M$ can be given the structure of a differentiable manifold,
since it can be identified with the quotient space $G/G_\gamma$,
where $G_\gamma$ is the closed subgroup of $G$ that fixes a given 
$\gamma$ in $\Xi$.

We will need a more precise description of the shortest closed geodesics and the Helgason spheres.
Let $p$ be some point in $M$. 
We write $M=G/K$ where $(G,K)$ is 
a symmetric pair and $K$ is the isotropy group of $G$ at $p$.
Consider  the corresponding Cartan decomposition
$\mathfrak{g}=\mathfrak{k}+\mathfrak{p}$
of the Lie algebra $\mathfrak{g}$ of $G$.
Then the tangent space $T_pM$ can be identified with the subspace
$\mathfrak{p}$ in the usual way. 
Let $\mathfrak{a}$ be a maximal 
Abelian subspace of $\mathfrak{p}$. Recall that a linear form $\alpha$ 
on $\mathfrak{a}$ is a restricted root if $\alpha \neq 0$ and
$\mathfrak{g}_{\alpha} \neq \{0\}$ hold, where the subspace
$\mathfrak{g}_{\alpha}$ is defined by
\begin{equation*} 
\mathfrak{g}_{\alpha}=\{\, X\in \mathfrak{g} \mid
(\ad H)^2(X)=-\alpha(H)^2 X\ \mathrm{for\ all}\ H\in \mathfrak{a} \,\}.
\end{equation*} 
Let $\mathcal{R}$ denote the set of roots.
Consider the subspaces
$\mathfrak{k}_{\alpha}=\mathfrak{g}_{\alpha}\cap \mathfrak{k}$
and
$\mathfrak{p}_{\alpha}=\mathfrak{g}_{\alpha}\cap \mathfrak{p}$
for $\alpha\in \mathcal{R}\cup \{0\}$. 
We obtain orthogonal decompositions
$\mathfrak{k}=\mathfrak{k}_0+\sum_{\alpha \in \mathcal{R}^{+}}
\mathfrak{k}_{\alpha}$
and
$\mathfrak{p}=\mathfrak{p}_0+\sum_{\alpha \in \mathcal{R}^{+}}
\mathfrak{p}_{\alpha}$
with respect to the Killing form $B$ of $\mathfrak{g}$,
where $\mathcal{R}^{+}$ is a set of  roots such that
$\mathcal{R}=\mathcal{R}^+\cup (-\mathcal{R}^+)$ is a disjoint 
decomposition of $\mathcal{R}$. 
Notice that $\mathfrak{a}=\mathfrak{p}_0$.

Consider now the exponential map of $M$ denoted by $\Exp$ and 
the lattice $\mathfrak{a}_K =\{\, X \!\in\! \mathfrak{a} \mid 
\exp(X) \!\in\! K \}$ in $\mathfrak{a}$.
It is clear that the flat torus $\Exp(\mathfrak{a})$
can be identified with the coset space
$\mathfrak{a}/\mathfrak{a}_K$ in the usual way.
Using the Killing form $B$, we can associate to each restricted
root $\alpha$ the vector $H_{\alpha} \!\in\! \mathfrak{a}$
for which $B(H_{\alpha},H)=\alpha(H)$
holds for any element $H$ of $\mathfrak{a}$.
Then by Theorem 8.5 in Chapter VII of \cite{He2}
the lattice $\mathfrak{a}_K$ is generated by the
vectors 
$$
X_{\alpha}=\frac{2\pi}{B(H_{\alpha},H_{\alpha})}\,H_{\alpha},
$$
$\alpha \in \mathcal{R}$.
It follows that the geodesic $\gamma_\alpha:[0,1]\to M$ defined by setting
$\gamma_\alpha(t)=\Exp(tX_\alpha)$ is closed and it is a shortest closed 
geodesic if $\alpha$ is a longest root.
Note that
these vectors 
$X_{\alpha},\ \alpha\in\mathcal{R}$, form
a dual root system of $\mathcal{R}$ in $\mathfrak{a}$.

The following proposition is proved in \cite{He2}, Chapter VII, \S11.

\begin{proposition} \label{Lietriple}
Let $\gamma_\delta$ be a closed geodesic corresponding to a 
longest root $\delta$ in $\mathcal{R}$. Then $\mathfrak{s}_\delta=\mathbb{R}H_\delta +\mathfrak{p}_\delta$ 
is a Lie triple system and $S=\Exp(\mathfrak{s}_\delta)$ is a Helgason sphere containing $\gamma_\delta$. 
\end{proposition}

Let $p$ be some point in $M$. The {\it midpoint locus of $p$} is the set $A_p$
of midpoints of shortest closed geodesics starting in $p$. 
It is an immediate consequence of Theorem \ref{helgspheres} (2) that $A_p$
coincides with the orbit $K(\gamma(1/2))$ where $K$ is the isotropy group of $G$
at $p$ and $\gamma:[0,1]\to M$ is some shortest closed geodesic starting at $p$.  The midpoint
loci are totally geodesic; see Corollary 11.13 in Chapter VII of \cite{He2}.
\medskip

The following proposition is crucial for the proof of Theorem \ref{support}. Notice that the
midpoint locus of a point on a sphere consists only of the antipodal point. 

\begin{proposition} \label{dim}
Let $p$ be some point in $M$. Then the dimension of the midpoint locus of $p$ is at least two
if $M$ is not a sphere.
\end{proposition}

\begin{proof}
It is well-known that the
proposition  is true
if $M$ is a compact symmetric space of rank one. We therefore assume that
$M$ is an irreducible compact symmetric space of rank greater than one.

Let $Y$ be an element of $\mathfrak{a}$ such that
$Y \in \frac{1}{2} \mathfrak{a}_K$ and $Y \notin \mathfrak{a}_K$ hold.
This means that $q=\Exp(Y)$ is an antipodal point of $p$
on the closed geodesic $\gamma(t)=\Exp(2tY)$, $t\in[0,1]$,
and $2\alpha(Y) \in \mathbb{Z}\pi$ holds for any
$\alpha\in\mathcal{R}$.
Then Theorem 11.14 in Chapter VII of \cite{He2} implies that 
the totally geodesic orbit $K(q)$ is isometric to
$\Exp(\sum_{\alpha\in\mathcal{R}(Y)} \mathfrak{p}_\alpha)$,
where 
$\mathcal{R}(Y) =\{\, \alpha \in \mathcal{R}^{+} \mid
\ 2\alpha(Y)\in(2\mathbb{Z}+1)\pi \,\}$.
Moreover, the Lie algebra of the isotropy group $K_q$ of $K$
at the point $q$ coincides with
$\mathfrak{k}_q=\mathfrak{k}_0 +
\sum_{\alpha\in\mathcal{R}^{+}\setminus\mathcal{R}(Y)}
\mathfrak{k}_{\alpha}$.
 
Select a longest restricted root $\delta$ and take the shortest
closed geodesic $\gamma_\delta:[0,1]\to M$ defined by 
$\gamma_\delta(t)=\Exp(tX_{\delta})$.
The irreducibility of $M$ implies that
the root system $\mathcal{R}$ is irreducible. Hence,
we can find two roots $\beta_1,\ \beta_2$ in $\mathcal{R}^{+}$ such that 
$X_{\beta_1}$ and $X_{\beta_2}$ are not perpendicular to $X_{\delta}$.
Then we obtain that
\begin{equation*} 
\beta_i(X_{\delta}) =
\pi \cdot \dfrac{2\,B(H_{\beta_i},H_\delta)}{B(H_\delta,H_\delta)}
=\pm \pi
\end{equation*} 
holds for $i=1,2$, where the second equals sign follows from the fact
that $\delta$ is one of the longest roots and therefore 
\,$\Vert H_\delta \Vert \geq \Vert H_{\beta_i} \Vert$ holds.
Therefore 
$\mathcal{R}(X_{\delta}/2)$ contains at least two restricted roots 
and the dimension of the midpoint locus
$A_p=K(\Exp(X_{\delta}/2))$ is at least two.
\end{proof}

\begin{remark} We can of course define a midpoint locus with respect to
closed geodesics $\gamma_\alpha$ belonging to roots $\alpha$ which are
not longest. In this case, the conclusion of Proposition \ref{dim} is
not necessarily true as the following examples show.

Let us consider a compact symmetric space $M$ whose restricted root system
is of type $B_2$; one can for example choose $M$ to be a complex quadric 
$Q^n$ in $P^{n+1}(\mathbb{C})$  with $n\ge 3$.
Select a basis $\alpha_1,\ \alpha_2$ of this root system $\mathcal{R}$,
where $\alpha_1$ is the longer of the two roots.
Then the other positive roots are $\beta=\alpha_1+\alpha_2$
and $\delta=\alpha_1+2\alpha_2$. Consider the dual root system in 
$\mathfrak{a}$ which is represented in Figure 1. This shows that
the arc length of the closed geodesic defined by $X_{\beta}$
is equal to $2\sqrt{2}\pi$. It follows from the proof above that
the orbit $K(\Exp(X_{\beta}/2))$ of the antipodal point
$\Exp(X_{\beta}/2)$ coincides with a single point.
\begin{figure}[tbh]
\begin{center}
\includegraphics[width=6.5cm]{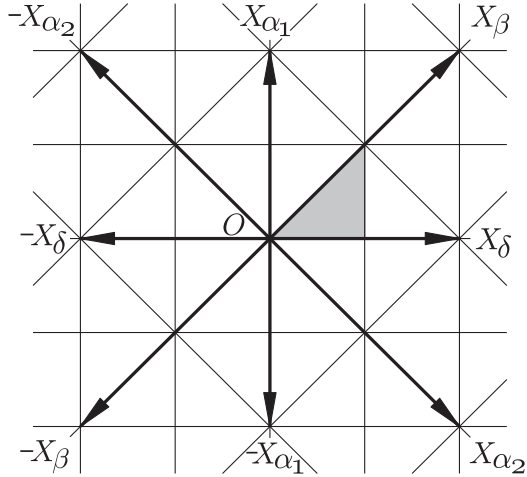}
\caption{The dual root vectors and the diagram of $M$ in $\mathfrak{a}$.}
\label{fig:fig1}
\end{center}
\end{figure}
\end{remark}

\section{Proof of Theorem \ref{theorem}} \label{sec:proof1}

Let $M$ be an irreducible simply connected compact symmetric space which is not a sphere.
Let $f$ be a smooth function on $M$ such that $\hat f=0$. Assume that $f$ does not vanish.
Then $f$ is not constant since otherwise $\hat f$ would clearly also be constant and
nonvanishing. Now it follows from Sard's Lemma that $f$ has a regular value $b$ in its image.
Hence $N=f^{-1}(b)$ is a hypersurface in $M$. 
Let $p$ be some point in $N$ and let $S$ be a Helgason sphere
containing $p$. 
Let $q$ be the antipodal point of $p$ in $S$. 
Since the great circles on $S$ are shortest closed geodesics in $\Xi$, it follows that $\widehat{f|S}=0$.
By Funk's Theorem \ref{Funk}, $f|S$ is an odd function. Hence $f(q)=-b$. 
Now let $g$ be an element of the isotropy group $G_p$ of $p$ and $A_q$ the midpoint locus of $q$.
We would like to prove the
inclusion
$$
g(A_q)=A_{g(q)}\subset N.
$$
First notice that
$g(q)$ is the antipodal point of $p$ in $g(S)$ since $g$ fixes $p$. The value of $f$ being $b$ in 
$p$ implies that it is equal to $-b$ in $g(q)$. This in turn implies that the value of $f$ on the 
midpoint locus $A_{g(q)}=g(A_q)$ is equal to $b$ since $A_{g(q)}$ consists of points that are 
antipodes of $g(q)$ in some Helgason sphere.
Hence we have proved that $g(A_q)\subset N$
for every $g$ in $G_p$ as we wanted.

Now let $v$ be a nonzero tangent vector in $T_pA_q$ which exists by Proposition \ref{dim}.
If $g$ is in the isotropy group $G_p$, then $dg_p(v)\in T_pg(A_q)\subset T_pN$
where $dg_p$ denotes the differential of $g$ at $p$. Hence the 
orbit $G_p(v)$ is contained in $T_pN$. Let $V$ be the linear subspace of $T_pN$ spanned
by the orbit $G_p(v)$. It is easy to see that $V$ is invariant under $G_p$. It is clear that
$V$ is  a proper subspace of $T_pM$ since it is contained in the hyperplane $T_pN$ in $T_pM$.
This means that the isotropy representation of $G_p$ is reducible which contradicts the irreducibility 
of the symmetric space $M$. This finishes the proof of Theorem \ref{theorem}. \qed

\begin{remark} \label{simple}
We give a very simple argument that can be used to prove Theorem \ref{theorem}
for the already known special case of simply connected projective spaces.
Note  that the simply connected projective spaces
with a symmetric metric
are nothing but the simply connected compact rank one spaces excluding the spheres. Let $M$ be  
such a space. A Helgason sphere in $M$ is nothing but a projective line. Let $f$ be a smooth
function on $M$ such that $\hat f$ vanishes. Then arguing as in the proof above we see that the
restrictions of $f$ to the projective lines are odd functions. Now let $p$ be some point in $M$. 
Let $P$ be a projective line through $p$. Let $q$ be the 
antipodal point of $p$ on $P$ and let $Q$ be a projective line that meets $P$ perpendicularly in
$q$. Let $r$ be the antipodal point of $q$ in $Q$ and let $R$ be the unique projective line that
connects $r$ and $p$. It is now an easy geometric exercise to show that $r$ and $p$ are antipodal 
on $R$.  Now we finish the proof as follows. Let $a$ be the value of $f$ in $p$. The value of
$f$ in $q$ is $-a$ since $p$ and $q$ are antipodal in $P$. The value of $f$ in $r$ is therefore
$a$ since $q$ and $r$ are antipodal in $Q$. Finally the value of $f$ in $p$ is $-a$ since 
$r$ and $p$ are antipodal in $R$. We have proved that the value of $f$ in $p$ is both $a$
and $-a$. Hence $a=0$. This proves that $f$ vanishes identically since $p$ is arbitrary. 

\end{remark}

\section{Proof of Theorem \ref{support}}

In Theorem \ref{support}, we are dealing with the simply connected symmetric
spaces of rank one which are not spheres.  These spaces are nothing but the
complex and quaternionic projective spaces and the octonion plane. The Helgason 
spheres in these spaces are precisely the projective lines.

The proof of Theorem \ref{support} will be very similar to the one of Theorem \ref{theorem} that
we just gave in Section \ref{sec:proof1}. Let $q$ be as in Theorem \ref{support}. The
assumptions in the theorem only give information on $q$ if there is a closed geodesic $\gamma$ in
$\Xi$ passing through $q$ and not meeting $\bar B_r(p)$. To be able to carry over the
ideas of the proof of Theorem \ref{theorem} we need a Helgason sphere $S$, which is here a
projective line, passing through $q$ and not meeting $\bar B_r(p)$.

\begin{proposition} \label{exprojline}
Let $M$ be a simply connected compact symmetric space of rank one which is not a sphere.
Let $\bar B_r(p)$ be a closed ball of radius $r>0$ around a point $p$ in $M$ where $r$ is
smaller than the injectivity radius $i(M)$ of $M$.  Let $q\not\in \bar B_r(p)$.
Then there is a projective line $S$ through $q$ that
does not meet the open ball $B_s(p)$ where $s=d(p,q)$.
In particular, $S$ does not meet $\bar B_r(p)$.
\end{proposition}

\begin{proof}  We have normalized the metric  on $M$ so  that the injectivity radius 
satisfies $i(M)=\pi$. The diameter $d(M)$ of $M$ then also satisfies $d(M)=\pi$
since $M$ has rank one and  is simply connected.

If $s=d(p,q)=\pi$, then $q$ is contained in the cut locus $C(p)$ of $p$
(which for simply connected rank one symmetric spaces coincides with the
midpoint locus $A_p$ of $p$).
 Note that   $C(p)$  is a projective line in $M$
if $M$ is a projective plane and a projective hyperplane in $M$ otherwise.
We can now choose a projective line $S$ in $C(p)$ passing through $q$ which then
clearly satisfies the claim in the
proposition since $d(p,C(p))=\pi$.

We now assume that $s=d(p,q)< \pi$. Let
$\sigma:[0,s]\to M$ be the shortest geodesic connection between $q$ and $p$.
Let $H$ be the hyperplane in $T_qM$ that is perpendicular to $\dot\sigma(0)$. 
We first show that there is a projective line $S$  passing through $q$ such that 
$T_qS\subset H$. We continue the geodesic $\sigma$ beyond $p$ until we reach
the first cut point $\hat p=\sigma(t_0)$ on $\sigma$.  Note that $t_0=\pi$.
We let $C(\hat p)$ denote the
cut locus of $\hat p$.
Note that $d(\hat p,C(\hat p))=\pi$.
As above we can choose $S$ as a projective line in $C(\hat p)$ that passes
through $q$. Clearly $S$ and the open ball $B_\pi(\hat p)$ do not meet.

We will prove that $S$ does not meet the open ball
$B_s(p)$. Assume that $\hat q\in S\cap B_s(p)$. There is a geodesic $\tau$ of
length less than $s$ from $p$ to $\hat q$. Then the concatenation $\hat\tau$
of $\sigma^{-1}|[s,\pi]$
and $\tau$ is a broken geodesic connecting $\hat p$ and $\hat q$. The length of
$\hat\tau$ is less than $\pi$ since $L(\sigma|[s,\pi])=\pi-s$. This is a contradiction
since $d(\hat p,S)=\pi$. Hence the claim of the proposition follows.
\end{proof}

The proof of Theorem \ref{support} will be very similar to the proof of Theorem \ref{theorem}
in Section \ref{sec:proof1}. 

\begin{proof}[Proof of Theorem \ref{support}]
Let $q$ be a point outside of $\bar B_r(p)$ and let $S$
a projective line passing through $q$ and not meeting $\bar B_r(p)$ which exists
by Proposition \ref{exprojline}.
We assume that $f(q)\not=0$. Set $f(q)=a$. The value of $f$ in the antipode of
$q$ on $S$ is $-a$. Hence
 there is a regular value $b$ in the interval $(-a,a)$
 and $N=f^{-1}(b)$ is a hypersurface in $M$ which meets $S$ in a point $\hat q$.
Let $\tilde q$ be the antipodal point of $\hat q$ in $S$. It follows that $f(\tilde q)=-b$.
Let $G_{\tilde q}$ be the isotropy group at $\tilde q$, and let $V$ be a neighborhood of 
the identity in $G_{\tilde q}$ with the property that $g(S)$ does not meet $\bar B_r(p)$
for any $g$ in $V$. Then $V\cdot \hat q$ is a neighborhood of $\hat q$ in the midpoint locus 
$A_{\tilde q}$
and it follows that $f|V\cdot \hat q$ is constant equal to $b$, i.e.,
$$
V\cdot \hat q\subset N.
$$
Now let $U$ be a neighborhood of the identity in the isotropy group $G_{\hat q}$ such that
$g(h(S))$ does not meet $\bar B_r(p)$ for any $g\in U$ and $h\in V$.
Then we can show with methods as in 
Section \ref{sec:proof1} that
$$
g(V\cdot \hat q)\subset N
$$
for every $g\in U$.

We choose a nonzero tangent vector $v$ in $T_{\hat q}A_{\tilde q}=T_{\hat q}(V\cdot \hat q)$ which exists
by Proposition \ref{dim}. Arguing as in Section \ref{sec:proof1}, we see that the 
neighborhood $U\cdot v$ of $v$ in $G_{\hat q}(v)$ is contained in $T_{\hat q}N$. Hence $G_{\hat q}(v)$
is contained in $T_{\hat q}N$. This 
is in contradiction to the irreducibility of $M$ and finishes the proof of Theorem
\ref{support}. 
\end{proof}

\section{Conjecture}

We conjecture that Theorems \ref{theorem} and \ref{support} can be generalized as follows:

\begin{conjecture} \label{conjecture}
Assume that   $M$ is not a sphere. 
Then there is a number $r_0(M)$ depending on $M$ in the
halfopen interval $(0,d(M)]$ with the following property: 
Let $\bar B_r(p)$ be a closed ball of radius $r< r_0(M)$ around a point $p$. Let $f$ be a function on $M$ such that 
$\hat f(\gamma)=0$ for every $\gamma$ in $\Xi$ which does not meet $\bar B_r(p)$. 
Then $f(q)=0$ for all $q\not\in \bar B_r(p)$.
\end{conjecture} 

\begin{remark} In Theorem \ref{support} we can choose $r_0(M)=d(M)$. In more general
symmetric spaces, we expect that $r_0(M)$ has to be chosen smaller than $d(M)$ and 
maybe also smaller than $i(M)$.
\end{remark}


\end{document}